\documentclass[reqno,12pt]{amsproc}
\newtheorem{theo}{Теорема}
\newtheorem{rem}{Замечание}

\newtheorem{lem}{Лемма}

\newtheorem{df}{Определение}

\newcommand\eps\varepsilon
\newcommand\ph\varphi
\newcommand\kap\varkappa

\newcommand\bs {\boldsymbol}

\usepackage{hyperref}
\usepackage{url}

\usepackage{enumerate}
\usepackage{graphicx}
\usepackage{float}
\usepackage{yfonts}
\usepackage{mathrsfs}
\usepackage{amsthm}
\usepackage{amsmath}
\usepackage{amscd}

\usepackage[cp1251]{inputenc}
\usepackage[english, russian]{babel}

\begin{document}\title
{Принцип Даламбера-Лагранжа и уравнения Лагранжа }
\author[Oleg Zubelevich]{Oleg Zubelevich\\
 Dept. of Theoretical mechanics,  \\
Mechanics and Mathematics Faculty,\\
M. V. Lomonosov Moscow State University\\
Russia, 119899, Moscow,  MGU \\
 }
\maketitle

\section{Описание задачи на физическом уровне строгости}
Пусть задана система материальных точек с массами
$m_1,\ldots,m_\nu.$
Через $\boldsymbol r_1,\ldots,\boldsymbol r_\nu$ обозначим радиус-векторы этих точек в некоторой инерциальной декартовой системе $OXYZ$:
$$\boldsymbol r_i=(X_i,Y_i,Z_i).$$
Силу, действующую на  $i-$ю точку, представим в виде суммы
$\boldsymbol F_i+\boldsymbol R_i,$ где
\begin{align}\boldsymbol F_i&=\boldsymbol F_i(t,\boldsymbol r_1,\ldots,\boldsymbol r_{\nu},\boldsymbol {\dot r}_1,\ldots,\boldsymbol {\dot r}_{\nu})=(F_i^x,F_i^y,F_i^z),\nonumber\\
 \boldsymbol R_i&=\boldsymbol R_i(t,\boldsymbol r_1,\ldots,\boldsymbol r_{\nu},\boldsymbol {\dot r}_1,\ldots,\boldsymbol {\dot r}_{\nu})=(R_i^x,R_i^y,R_i^z).\nonumber
 \end{align}
Силы $\boldsymbol F_i$ будем cчитать заданными и называть активными силами. 

Силы $\boldsymbol R_i$ принуждают систему материальных точек двигаться в соответствии с уравнениями
$$\Phi^k(t,\boldsymbol r_1,\ldots,\boldsymbol r_{\nu},\boldsymbol {\dot r}_1,\ldots,\boldsymbol {\dot r}_{\nu})=0,\quad k=1,\ldots, n,$$ которые называются уравнениями связей или связями. Силы $\boldsymbol R_i$ называются силами реакции связей или реакциями. Реакции  подлежат определению в ходе решения задачи. Если связи отсутствуют, то силы реакций полагают равными нулю.

Запишем второй закон Ньютона для каждой точки:
\begin{equation}\label{dfgjjj}m_i\boldsymbol {\ddot r}_i=\boldsymbol F_i+\boldsymbol R_i,\quad i=1,\ldots, {\nu}.\end{equation}
Чтобы перейти к аккуратным формулировкам, введем обозначения:
\begin{align}G&=\mathrm{diag}\,(m_1,m_1,m_1,\ldots,m_{\nu},m_{\nu},m_{\nu});\quad x=(X_1,Y_1,Z_1,\ldots,X_{\nu},Y_{\nu},Z_{\nu})^T;\nonumber\\
f(t,x,\dot x)&=(F_1^x,F_1^y,F_1^z,\ldots ,F_{\nu}^x,F_{\nu}^y,F_{\nu}^z),\quad N(t,x,\dot x)=(R_1^x,R_1^y,R_1^z,\ldots ,R_{\nu}^x,R_{\nu}^y,R_{\nu}^z),\nonumber\\
\ph^k(t,x,\dot x)&=\Phi^k.\nonumber
\end{align}

\section{Принцип Даламбера-Лагранжа }
Вектор-строка $f$ и вектор-столбец $\ph=(\ph^1,\ldots,\ph^n)^T$ предполагаются гладкими функциями\footnote{Далее мы считаем все функции настолько гладкими, насколько это необходимо, чтобы формулы имели смысл.} на прямом произведении $I\times D\times \mathbb{R}^m.$ Здесь $D\subset\mathbb{R}^m$  -- область, $$ m=3\nu,\quad x=(x^1,\ldots, x^m)^T\in D,\quad t\in I=(t_1,t_2),\quad \dot x\in\mathbb{R}^m.$$
 Система  (\ref{dfgjjj}) записывается следующим образом:
\begin{equation}\label{seg44}
G\ddot x=f^T(t,x,\dot x)+N^T(t,x,\dot x),\end{equation}
 а связи принимают вид:
\begin{equation}\label{xfg55}
\ph(t,x,\dot x)=0.\end{equation}
Область $D$ называется конфигурационным пространством; область $D\times\mathbb{R}^m\ni (x,\dot x)$ называется фазовым пространством; область $I\times D\times\mathbb{R}^m\ni (t,x,\dot x)$ называется расширенным фазовым пространством системы (\ref{seg44}).

Представим систему (\ref{seg44})  в ковариантной форме:
\begin{equation}\label{azzxxx}[\mathscr T]=f+N,\quad[\mathscr T]= \frac{d}{dt}\frac{\partial\mathscr{T}}{\partial \dot x}-\frac{\partial\mathscr{T}}{\partial  x}, \end{equation}
где
$$\mathscr T=\frac{1}{2}\sum_{i=1}^\nu m_i|\boldsymbol {\dot r}_i|^2=\frac{1}{2}\dot x^TG\dot x$$
-- кинетическая энергия системы.
\begin{rem}Представление системы (\ref{seg44})  в форме (\ref{azzxxx})
играет ключевую роль в дальнейшем изложении (см. раздел \ref{cbbbbb}). Поэтому  мы будем придерживаться этой формы, хотя на данном этапе это может показаться странным.\end{rem}

Если существует функция $\mathscr V(t,x,\dot x)$ такая, что
$f=[\mathscr{V}],$
то активные силы называются обобщенно потенциальными.
Функция $\mathscr V(t,x,\dot x)$ называется обобщенным потенциалом.
Если обобщенный потенциал не зависит от $\dot x$, то он называется потенциалом или потенциальной энергией.

Предположим, что связи невырожденны:
$$ \mathrm{rang}\,\ph_{\dot x}=n<m,\quad\ph_{\dot x}=\begin{pmatrix}
\frac{\partial \ph^1}{\partial \dot x^1}  & \cdots & \frac{\partial \ph^1}{\partial \dot x^m} \\
\vdots    & \ddots & \vdots  \\
\frac{\partial \ph^n}{\partial \dot x^1} & \cdots & \frac{\partial \ph^n}{\partial \dot x^m}
\end{pmatrix}.$$
Множество решений уравнения (\ref{xfg55})
$$S=\{(t,x,\dot x)\mid\ph(t,x,\dot x)=0\},\quad \dim S=2m+1-n,$$ которое мы считаем непустым,  является гладким подмногообразием расширенного фазового пространства.

\begin{theo}[Принцип освобождаемости от связей]\label{as11}
Существует, и притом единственная, гладкая на $I\times D\times \mathbb{R}^m$ вектор-строка $N=N(t,x,\dot x)$ такая, что

1) функции $\ph$ являются первыми интегралами системы (\ref{azzxxx});

2) для всякого вектор-столбца $\xi\in\mathbb{R}^m$ такого, что
\begin{equation}\label{sdg-0}\ph_{\dot x}(t,x,\dot x)\xi=0\end{equation}
выполнено равенство
\begin{equation}\label{sfg43}N(t,x,\dot x)\xi =0.\end{equation}
\end{theo}\begin{df}Линейное пространство векторов $\xi\in\mathbb{R}^m$, удовлетворяющих (\ref{sdg-0}), т.е. $\ker\ph_{\dot x}(t,x,\dot x)$, называется пространством виртуальных перемещений. Это пространство отнесено к значениям переменных $(t,x,\dot x)$.

Размерность пространства виртуальных перемещений называется числом степеней свободы системы
(\ref{azzxxx}), (\ref{xfg55}).
\end{df} Очевидно, число степеней свободы равно $m-n=:r$.

\begin{df} Если силы $N$ выбраны в соответствии с теоремой \ref{as11}, то система (\ref{azzxxx}), (\ref{xfg55}) называется системой с идеальными связями. Связи (\ref{xfg55}) называются идеальными, а силы $N$ называются реакциями идеальных связей.\end{df}

Ниже обсуждаются только системы с идеальными связями.

 По теореме \ref{as11}, многообразие $S$ является инвариантным многообразием системы (\ref{azzxxx}).

\subsubsection*{ Доказательство теоремы \ref{as11}.} \begin{lem}\label{dfh6ll} Пусть $X,Y,Z$ -- векторные пространства и
$$A:X\to Y,\quad B:X\to Z$$ -- линейные операторы.
Тогда если $\ker A\subset\ker B,$ то cуществует  оператор $\Lambda:Y\to Z$ такой, что $B=\Lambda A$.
\end{lem}
Отсюда и из формул (\ref{sdg-0}), (\ref{sfg43}) следует, 
что найдется вектор-строка $\Lambda=(\lambda_1,\ldots,\lambda_n)(t,x,\dot x)$ такая, что 
\begin{equation}\label{dfg54hh}N=\Lambda\ph_{\dot x}.\end{equation}
Продифференцируем $\ph$ в силу системы (\ref{azzxxx}) и приравняем результат к нулю:
$$\ph_t+\ph_x\dot x+\ph_{\dot x}G^{-1}(f^T+N^T)=0.$$
Подставляя сюда $N$ из формулы  (\ref{dfg54hh}), находим
\begin{equation}\label{xdfgaxv}\Lambda^T=-\big(\ph_{\dot x}G^{-1}\ph_{\dot x}^T\big)^{-1}(\ph_t+\ph_x\dot x+\ph_{\dot x}G^{-1}f^T).\end{equation}
Здесь мы используем следующий факт из линейной алгебры.
\begin{lem}\label{xfg6zzzz6}Если квадратная $m\times m-$матрица $P$ симметрична и положительно определена, а матрица $B$ состоит из $m$ столбцов и $n$ строк ($m>n$), причем $\mathrm{rang}\,B=n,$ то матрица $BPB^T$ симметрична и положительно определена.\end{lem}
Теорема доказана.

Уравнение (\ref{azzxxx}) с подставленной в него реакцией  (\ref {dfg54hh}) называется уравнением Лагранжа со множителями или уравнением Лагранжа первого рода.

\begin{theo}\label{dfg44}
Пусть функция $x(t)$ удовлетворяет (\ref {xfg55}) и 
\begin{align}\ph_{\dot x}&\big(t,x(t),\dot x(t)\big)\xi=0\Longrightarrow \nonumber\\
&([\mathscr T]-f)\xi=0.\label{xsdfg23ee}\end{align}
Тогда $x(t)$ удовлетворяет (\ref{azzxxx}) с реакциями (\ref{dfg54hh}), (\ref{xdfgaxv}).
\end{theo}Уравнение (\ref{xsdfg23ee}) с виртуальными перемещениями $\xi$ называется общим уравнением динамики.

Если последовательно подставить в (\ref{xsdfg23ee}) векторы базиса пространства виртуальных перемещений, то получится система из $r$ скалярных дифференциальных уравнений, каждое из которых содержит компоненты вектора $\ddot x$. Вместе с уравнениями связей получится система из $m$ дифференциальных уравнений порядка $2r+n=2m-n=\dim S-1$.

{\it Доказательство теоремы \ref{dfg44}.} Действительно, в силу леммы \ref{dfh6ll} существует вектор-строка
$\tilde\Lambda(t)$ такая, что 
\begin{equation}\label{dfh6pp}
[\mathscr T]-f=\ddot x^T(t)G-f\big(t,x(t),\dot x(t)\big)=\tilde\Lambda(t)\ph_{\dot x}\big(t,x(t),\dot x(t)\big).\end{equation}
Продифференцируем уравнение $\ph\big(t,x(t),\dot x(t)\big)=0:$
$$\ph_t\big(t,x(t),\dot x(t)\big)+\ph_x\big(t,x(t),\dot x(t)\big)\dot x(t)+\ph_{\dot x}\big(t,x(t),\dot x(t)\big)\ddot x(t)=0.$$
Подставляя сюда $\ddot x$ из (\ref{dfh6pp}), убеждаемся, что 
$\tilde \Lambda(t)=\Lambda\big(t,x(t),\dot x(t)\big).$

Теорема доказана.

Следующая теорема очевидна.
\begin{theo}\label{dfg4xde4}
Пусть $x(t)$ удовлетворяет (\ref{azzxxx}) с реакциями (\ref{dfg54hh}), (\ref{xdfgaxv}).
Тогда верна импликация $$\ph_{\dot x}\big(t,x(t),\dot x(t)\big)\xi=0\Longrightarrow
([\mathscr T]-f)\xi=0.$$\end{theo}

Теоремы \ref{as11}, \ref{dfg44}, \ref{dfg4xde4} составляют принцип Даламбера-Лагранжа.

\section{О корректности определения  реакций и виртуальных перемещений}\label{sdf55}
 В определения реакций $N$ и пространства виртуальных перемещений входят функции $\ph$.
В этом разделе мы покажем, что ни реакции, ни пространство виртуальных перемещений  не зависят от выбора функций $\ph$, задающих многообразие $S$. Другими словами, реакции и виртуальные перемещения зависят от геометрии многообразия $S$, а не от конкретного способа его аналитического задания.

Введем вектор гладких функций, определенных на $I\times D\times\mathbb{R}^m\times\mathbb{R}^n:$
$$U=(U^1,\ldots,U^n)^T(t,x,\dot x,z),\quad z\in\mathbb{R}^n,$$
такой, что
$$U(t,x,\dot x,z)=0\Leftrightarrow z=0,\quad \det U_z(t,x,\dot x,0)\ne 0,\quad \forall (t,x,\dot x)\in I\times D\times\mathbb{R}^m.$$
Положим $\psi(t,x,\dot x):=U\big(t,x,\dot x,\ph(t,x,\dot x)\big).$ Тогда $S=\{\ph=0\}=\{\psi=0\}$
 и верны формулы
$$\psi_t=U_t(t,x,\dot x,\ph)+U_z(t,x,\dot x,\ph)\ph_t,\quad U_t(t,x,\dot x,\ph)\mid_{\ph=0}=0.$$
Аналогично
$\psi_x\mid_S=U_z\ph_x\mid_S,\quad \psi_{\dot x}\mid_S=U_z\ph_{\dot x}\mid_S.$

В частности $$\mathrm{rang}\,\big(\psi_{\dot x}\mid_S\big)=n, \quad \ker \big(\psi_{\dot x}\mid_S\big)=\ker\big( \ph_{\dot x}\mid_S\big).$$

Вычислим реакции, отвечающие разным способам задания многообразия $S$:
\begin{align}
N^T&=-\ph^T_{\dot x}\big(\ph_{\dot x}G^{-1}\ph_{\dot x}^T\big)^{-1}(\ph_t+\ph_x\dot x+\ph_{\dot x}G^{-1}f^T)\nonumber\\
\tilde N^T&=-\psi^T_{\dot x}\big(\psi_{\dot x}G^{-1}\psi_{\dot x}^T\big)^{-1}(\psi_t+\psi_x\dot x+\psi_{\dot x}G^{-1}f^T).\nonumber
\end{align}
Отсюда
$
N^T\mid_S=\tilde N^T\mid_S.$
\section{Голономные и неголономные связи}
\begin{df}Связи (\ref{xfg55}) называются голономными, если найдутся функции
$$g=(g^1,\ldots, g^n)^T(t,x),\quad \mathrm{rang}\, g_x=n$$ такие, что 
$$S=\Big\{(t,x,\dot x)\mid\frac{d}{dt}g(t,x) =0\Big\},\quad \frac{d}{dt}g(t,x)=g_t(t,x)+g_x(t,x)\dot x.$$

Связи, не являющиеся голономными, называются неголономными.
\end{df}
Вопрос о голономности связей эквивалентен вопросу об интегрируемости дифференциальной системы на многообразии \cite{stern}.

Уравнения $g(t,x)=0$ называются геометрическими связями, а уравнения (\ref{xfg55}) -- дифференциальными.

Далее, до конца раздела \ref{cbbbbb} мы будем считать связи голономными и писать
\begin{equation}\label{xdfg11gy}\ph(t,x,\dot x)=g_t(t,x)+g_x(t,x)\dot x=0.\end{equation}
В частности,
\begin{equation}\label{xf11q}\ph_{\dot x}=g_x,\quad\ker \ph_{\dot x}(t,x,\dot x)=\ker g_x(t,x).\end{equation}

\section{Ковариантность вариационной производной}
Предположим, что многообразие
$$\Sigma_t=\{x\in D\mid g(t,x)=0\},\quad \dim\Sigma_t=r$$ является образом гладкого многообразия $Y,\quad \dim Y=r$ при вложении\footnote{Вложением мы называем   гладкое инъективное отображение
$u(t,\cdot):Y\to\mathbb{R}^m$ такое, что для каждой  точки $y\in Y$ и ее образа $x=u(t,y)\in \mathbb{R}^m$ существуют окрестности $U_y\subset Y$ и $U_x\subset\mathbb{R}^m$ с локальными координатами $(y^1,\ldots,y^r)$ и $(x^1,\ldots,x^m)$ соответственно, в которых $u$ имеет вид: $x^i=y^i,\quad i=1,\ldots, r,\quad x^s=0,\quad s=r+1,\ldots,m,$ причем $U_x\cap \Sigma_t=\{x^s=0,\quad s=r+1,\ldots,m\}.$
 Указанные локальные координаты гладко зависят от $t\in I$.}
$$u(t,\cdot):Y\to D,\quad \Sigma_t=u(t,Y).$$
Через $y=(y^1,\ldots,y^r)^T$ обозначим локальные координаты в $Y$;
\begin{equation}\label{fg11}g(t,u(t,y))=0,\quad u=(u^1,\ldots,u^m)^T, \quad x=u(t,y).\end{equation}
Линейное отображение
$$u_y(t,y):T_yY\to T_{u(t,y)}\Sigma_t$$ является изоморфизмом, причем касательное пространство
$T_{u(t,y)}\Sigma_t$ совпадает с пространством виртуальных перемещений:
$$T_{u(t,y)}\Sigma_t=\ker \big(\ph_{\dot x}\mid_{x=u(t,y)}\big).$$
Это ясно, если продифференцировать формулу (\ref{fg11}) по $y$ и использовать (\ref{xf11q}):
$$g_x(t,u(t,y))u_y(t,y)=0.$$
Рассмотрим гладкую функцию $\mathscr{F}(t,x,\dot x)$, определенную на расширенном фазовом пространстве.
\begin{df}
Вариационной производной $\mathscr{F}$ называется следующий дифференциальный оператор
$$[\mathscr{F}]=([\mathscr{F}]_{1},\ldots ,[\mathscr{F}]_{m}),\quad [\mathscr{F}]_{k}=\frac{d}{dt}\frac{\partial \mathscr{F}}{\partial\dot x^k}-\frac{\partial \mathscr{F}}{\partial x^k}.$$
\end{df}
Очевидно, вариационная производная является линейной операцией и
для любой функции $w=w(t,x),\quad w:I\times D\to \mathbb{R}$
выполнено тождество 
$$\Big[\frac{dw}{dt}\Big]=0.$$

\begin{theo}\label{dfggg456}
Верно равенство
$$[\mathscr{F}]\Big|_{x=u(t,y)}u_y(t,y)=[F],$$
где функция $F:I\times TY\to\mathbb{R}$ ($TY$ -- касательное расслоение) определяется формулой
\begin{equation}\label{dff}F(t,y,\dot y)=\mathscr{F}\Big|_{x=u(t,y)}=\mathscr{F}\big(t,u(t,y),u_t(t,y)+u_y(t,y)\dot y\big).\end{equation}
\end{theo}

{\it Доказательство теоремы \ref{dfggg456}.} Проведем вычисления в координатах. По формуле (\ref{dff}) имеем:
$$\frac{\partial F}{\partial y^i}=\frac{\partial \mathscr{F}}{\partial x^p}\frac{\partial u^p}{\partial y^i}+
\frac{\partial \mathscr{F}}{\partial \dot x^p}\Big(\frac{\partial^2 u^p}{\partial y^l\partial y^i}\dot y^l+\frac{\partial^2u^p}{\partial t\partial y^i}\Big);$$
и
$$
\frac{d}{dt}\frac{\partial F}{\partial \dot y^i}=\frac{d}{dt}\Big(\frac{\partial \mathscr{F}}{\partial \dot x^p}
\frac{\partial u^p}{\partial y^i}\Big)=\frac{\partial u^p}{\partial y^i}\Big(\frac{d}{dt}\frac{\partial \mathscr{F}}{\partial \dot x^p}\Big)+\frac{\partial \mathscr{F}}{\partial \dot x^p}\Big(\frac{\partial^2 u^p}{\partial y^l\partial y^i}\dot y^l+\frac{\partial^2u^p}{\partial t\partial y^i}\Big).$$
Вычитая из одного равенства другое, получаем требуемый результат.
Теорема доказана.

\section{Уравнения Лагранжа }\label{cbbbbb}
Теорема \ref{dfggg456} означает, что  вариационная производная ведет себя как ковекторное поле относительно  отображения $u$.
Говоря  неформально, вывод уравнений Лагранжа состоит в том, что с помощью отображения $u$ мы осуществляем <<pullback>> общего уравнения динамики с области $D$ на  многообразие $Y$.

Перейдем к точным формулировкам.

\begin{df}Обобщенными силами называются компоненты вектор-строки
\begin{equation}\label{sfg40}Q(t,y,\dot y):=f\big(t,u(t,y),u_t(t,y)+u_y(t,y)\dot y\big) u_y(t,y).\end{equation}
\end{df}
\begin{theo}\label{xsfgtyuio}
Пусть функция $y(t)$ является решением уравнений Лагранжа (второго рода)
\begin{equation}\label{sdfgqas}[T]=Q,\end{equation}
где кинетическая энергия $T=T(t,y,\dot y)$ определена формулой:
\begin{equation}\label{xdfsgggtyu}
T=\mathscr T\mid_{x=u(t,y)}=
\frac{1}{2}\big(u_t+u_y\dot y\big)^TG\big(u_t+u_y\dot y\big).
\end{equation}

Тогда функция $x(t)=u(t,y(t))$ удовлетворяет общему уравнению динамики (\ref{xsdfg23ee}) со связями (\ref{xdfg11gy}).\end{theo}
{\it Доказательство теоремы \ref{xsfgtyuio}.} По теореме \ref{dfggg456} имеем
\begin{equation}\label{xdfg6yiii}([\mathscr{T}]-f)\Big|_{x=u(t,y(t))}u_y(t,y(t))=([T]-Q)\Big|_{y=y(t)}.\end{equation}
Правая часть этого равенства равна нулю по условию. Как уже отмечалось, образ оператора $u_y$ совпадает с пространством виртуальных перемещений, поэтому выполнено общее уравнение динамики.

Теорема доказана.
\begin{df}Многообразие $Y$ называется конфигурационным многообразием, а локальные координаты  $y$ называются обобщенными координатами системы (\ref{sdfgqas}).

Многообразия $TY,\quad I\times TY$ называются фазовым пространством и расширенным фазовым пространством соответственно.

Компоненты вектора $\dot y$ называются обобщенными скоростями.
\end{df}

\begin{theo}\label{dfb00dd}
Пусть функция $x(t)$ удовлетворяет общему уравнению динамики (\ref{xsdfg23ee}) со связями (\ref{xdfg11gy}), и при некотором $t_0\in I$ верно равенство
$g(t_0,x(t_0))=0.$

Тогда существует, и притом единственное, решение $y(t)$ уравнений (\ref{sdfgqas})  такое, что $$x(t)=u(t,y(t)).$$\end{theo}
{\it Доказательство теоремы \ref{dfb00dd}}.
Проинтегрируем равенство $$\frac{d}{dt}g(t,x(t))=0$$ в пределах от $t_0$ до $t:\quad g(t,x(t))-g(t_0,x(t_0))=0.$
Таким образом, $x(t)\in \Sigma_t,\quad t\in I$.

Для любого $t$ уравнение \begin{equation}\label{dh0009}x(t)=u(t,y)\end{equation} однозначно разрешимо относительно $y$, поскольку $u(t,\cdot):Y\to\Sigma_t$ -- диффеоморфизм. Обозначим это решение за $y(t),\quad x(t)=u(t,y(t))$.

Покажем, что $y(t)$ -- гладкая функция. Действительно,  (\ref{dh0009}) представляет собой систему из $m$ уравнений с $r-$мерным неизвестным вектором $y$.  Выделим из этой системы $r$ независимых уравнений и применим теорему о неявной функции в части гладкости решения по параметру.

Левая часть  равенства (\ref{xdfg6yiii}) равна нулю по условию.
Теорема доказана.

\begin{theo}\label{asdf4eew}Верна формула
$T=T_2+T_1+T_0,$ где $$T_2=\frac{1}{2}\dot y^Tu_y^TGu_y\dot y$$ -- положительно определенная квадратичная форма обобщенных скоростей; $T_1=u_t^TGu_y\dot y$ -- линейная форма, $T_0=u_t^TGu_t/2$.\end{theo}
Это утверждение  следует непосредственно из формулы (\ref{xdfsgggtyu}) и леммы \ref{xfg6zzzz6}.

Из теоремы (\ref{asdf4eew}) следует, что система (\ref{sdfgqas}) представима в нормальной форме:
$$\ddot y=a(t,y,\dot y),$$ и потому для нее справедлива теорема существования и единственности Коши.

Если силы $f$ обобщенно потенциальны, то и силы $Q$ также обобщенно потенциальны -- в силу ковариантности операции $[\cdot]$. Причем
$Q=[V],\quad V(t,y,\dot y)=\mathscr V\mid_{x=u(t,y)},$ и уравнения (\ref{sdfgqas}) превращаются в уравнения Лагранжа:
$$[L]=0,\quad L=T-V.$$Функция $L=L(t,y,\dot y)$ называется функцией Лагранжа или лагранжианом.

\section{О классических обозначениях}
Традиционно в механике вектор виртуальных перемещений $\xi$ записывается  так:
$$\xi=(\delta \boldsymbol r_1,\ldots,\delta\boldsymbol  r_\nu)^T.$$
Соответственно,  уравнения (\ref{sdg-0}) и (\ref{sfg43}) имеют вид
$$\sum_{i=1}^\nu\Big(\frac{\partial \Phi^k}{\partial\boldsymbol{\dot r}_i},\delta \boldsymbol r_i\Big)=0,\quad\sum_{i=1}^\nu(\boldsymbol R_i,\delta \boldsymbol r_i)=0,\quad k=1,\ldots,n.$$

Общее уравнение динамики (\ref{xsdfg23ee}) имеет вид
$$\sum_{i=1}^\nu(m_i\boldsymbol{\ddot r}_i-\boldsymbol F_i,\delta\boldsymbol r_i)=0.$$

Вложение $u$ вводят так:
$$u(t,\cdot):y\mapsto (\boldsymbol r_1(t,y),\ldots,\boldsymbol r_\nu(t,y)).$$
Соответственно, обобщенные силы (\ref{sfg40})
вычисляются по формуле:
$$Q_j=\sum_{i=1}^\nu\Big(\bs F_i,\frac {\partial\bs r_i}{\partial  y^j}\Big),\quad Q=(Q_1,\ldots,Q_r).$$

\end{document}